\documentstyle[12pt,amscd]{amsart}
\newtheorem{Theorem}{Theorem}[section] \newtheorem{Lemma}[Theorem]{Lemma}
\newtheorem{Corollary}[Theorem]{Corollary} \newtheorem{Proposition}[Theorem]{Proposition}

\def\height{\operatorname{ht}}\def\sup{\operatorname{sup}} \def\In{\operatorname{in}} \def\reg{\operatorname{reg}}  \def\sk{\smallskip\par} \def\To{\longrightarrow} 

\begin{document}
\title{Evaluations of initial ideals  and \smallskip\\ Castelnuovo-Mumford regularity}
\author{Ng\^o Vi\^et Trung}
\address{Institute of Mathematics, Box 631, B\`o H\^o, Hanoi, Vietnam}
\email{nvtrung@@hn.vnn.vn}
\begin{abstract} This paper characterizes the Castelnuovo-Mumford regularity by evaluating the initial ideal with respect to the reverse lexicographic order. \end{abstract}
\thanks{The author is partially supported by the National Basic Research Program of Vietnam}
\maketitle

\section{Introduction} 

Let $S = k[x_1,\ldots,x_n]$ be a polynomial ring over a field $k$ of arbitrary characteristic.
Let $I \subset S$ be an arbitrary  homogeneous ideal and 
$$0 \To F_p \To \cdots \To F_1 \To F_0 \To S/I \To 0.$$
a graded minimal free resolution of $S/I$. Write $b_i$ for the maximum degree of the generators of $F_i$. The {\it Castelnuovo-Mumford regularity} 
$$\reg(S/I) := \max\{b_i-i|\ i = 0,\ldots,p\}$$
is a measure for the complexity of $I$ in computational problems [EG], [BM], [V]. One can use Buchsberger's syzygy algorithm to compute $\reg(S/I)$. However, such a computation is often very big. Theoretically,  if char($k$) = 0, $\reg(S/I)$ is equal to the largest degree of the generators of the generic initial ideal of $I$ with respect to the reverse lexicographic order [BS]. But it is difficult to know when an initial ideal is generic. Therefore, it would be of interest to have other methods for the computation of $\reg(S/I)$. \sk

The aim of this paper is to present a simple method for the computation of $\reg(S/I)$ which is based only on evaluations of $\In(I)$, where $\In(I)$ denotes the initial ideal of $I$ with respect to the reverse lexicographic order. We are inspired by a recent paper of Bermejo and Gimenez [BG] which gives such a method for the computation of the Castelnuovo-Mumford regularity of projective curves. \sk

Let $d = \dim S/I$. For $i = 0,\ldots,d$ put $S_i = k[x_1,\ldots,x_{n-i}]$. Let $J_i$ be the ideal of $S_i$ obtained from $\In(I)$ by the evaluation $x_{n-i+1} = \cdots = x_n = 0$. Let $\tilde J_i$ denote the ideal of $S_i$ obtained from $J_i$ by the evaluation $x_{n-i} = 1$. These ideals can be easily computed from the generators of $\In(I)$. In fact, if $\In(I) = (f_1,\ldots,f_s)$, where $f_1,\ldots,f_s$ are monomials in $S$, then $J_i$ is generated by the monomials $f_j$ not divided by any of the variables $x_{n-i+1},\ldots,x_n$ and $\tilde J_i$ by those monomials obtained from the latter by setting $x_{n-i} = 1$. Put 
$$c_i(I)  :=  \sup\{r|\ (\tilde J_i/J_i)_r \neq 0\},$$
with $c_i(I) = -\infty$ if $\tilde J_i = J_i$ and
$$r(I)  :=  \sup\{r|\ (S_d/J_d)_r \neq 0\}.$$

We can express $\reg(S/I)$ in terms of these numbers as follows. Assume that $c_i(I) < \infty$ for $i = 0,\ldots,d-1$. Then
$$\reg(S/I) = \max\{c_0(I),\ldots,c_{d-1}(I),r(I)\}.$$
The assumption $c_i(I) < \infty$ for $i = 0,\ldots,d-1$ is satisfied for a sufficiently general choice of the variables. If $I$ is the defining saturated ideal of a projective (not necessarily reduced) curve, this assumption is automatically satisfied if $k[x_{n-1},x_n]$ is a Noether normalization of $S/I$. In this case, $c_0(I) =-\infty$ and $\reg(S/I) = \max\{c_1(I),r(I)\}.$ From this formula we can easily deduce the results of Bermejo and Gimenez. \sk

Similarly we can compute the {\it partial regularities} 
$\ell\text{--}\!\reg(S/I) := \max\{b_i-i|\ i \ge \ell\}$, $\ell > 0$, which were recently introduced by Bayer, Charalambous and Popescu [BCP] (see also Aramova and Herzog [AH]).  These regularities can be defined in terms of local cohomology. Let $\frak m$ denote the maximal homogeneous ideal of $S$. Let $H_{\frak m}^i(S/I)$ denote the $i$th local cohomology module of $S/I$ with respect to $\frak m$  and set $a_i(S/I) = \max\{r|\ H_{\frak m}^i(S/I)_r \neq 0\}$ with $a_i(S/I) = -\infty$ if $H_{\frak m}^i(S/I) = 0$. For $t \ge 0$ we define
$\reg_t(S/I) := \max\{a_i(S/I)+i|\ i = 0,\ldots,t\}.$ Then 
$\reg_t(S/I) = (n-t)\text{--}\!\reg(S/I)$ [T2]. Under the assumption $c_i(I) < \infty$ for $i = 0,\ldots,t$ we obtain the following formula:
$$\reg_t(S/I) = \max\{c_i(I)|\ i = 0,\ldots,t\}.$$ 

The numbers $c_i(I)$ also allow us to determine the place at which $\reg(S/I)$ is attained in the minimal free resolution of $S/I$. In fact, $\reg(S/I) = b_t-t$ if $c_t(I) = \max\{c_i(I)|\ i = 0,\ldots,d\}$. Moreover, $r(I)$ can be used to estimate the reduction number of $S/I$ which is another measure for the complexity of $I$ [V]. \sk

It turns out that the numbers $c_i(I)$ and $r(I)$ can be described combinatorially in terms of the lattice vectors of the generators of $\In(I)$ (see Propositions 4.1-4.3 for details). These descriptions together with the above formulae give an effective method for the computation of $\reg(S/I)$ and $\reg_t(S/I)$. From this we can derive the estimation 
$$\reg_t(S/I) \le \max\{\deg g_i-n+i|\ i=0,\ldots,t\},$$
where $g_i$ is the least common multiple of the minimal generators of $\In(I)$ which are not divided by any of the variables $x_{n-i+1},\ldots,x_n$. \sk

This paper is organized as follows. In Section 2 we prepare some facts on the Castelnuovo-Mumford regularity. In Section 3 we prove the above formulae for $\reg(S/I)$ and $\reg_t(S/I)$. The combinatorial descriptions of $c_i(I)$ and $r(I)$ are given in Section 4. Section 5 deals with the case of projective curves. \sk

\noindent{\it Acknowledgement}. The author would like to thank M.~Morales for raising his interest in the paper of Bermejo and Gimenez [BG] and L.T.~Hoa for useful suggestions. \sk

\section{Filter-regular sequence of linear forms}

We shall keep the notations of the preceding section. Let ${\bf z} = z_1,\ldots,z_{t+1}$ be a sequence of homogeneous elements of $S$, $t \ge 0$. We call  $\bf z$ a {\it filter-regular sequence} for $S/I$ if $z_{i+1} \not\in {\frak p}$ for any associated prime ${\frak p} \neq {\frak m}$ of $(I,z_1,\ldots,z_i)$, $i = 0,\ldots,t$. This notion was introduced in order to characterize generalized Cohen-Macaulay rings [STC]. Recall that $S/I$ is a generalized Cohen-Macaulay ring if and only if $I$ is equidimensional and $(R/I)_{\frak p}$ is a Cohen-Macaulay ring for every prime ideal ${\frak p} \neq {\frak m}$. This condition is satisfied if $I$ is the defining ideal of a projective curve. We call $\bf z$ a homogeneous system of parameters for $S/I$ if $t+1 = d$ and $(I,z_1,\ldots,z_d)$ is an $\frak m$-primary ideal. It is known that every homogeneous system of parameters for $S/I$ is a filter-regular sequence if $S/I$ is a generalized Cohen-Macaulay ring. In general, a homogeneous system of parameters needs not to be a filter-regular sequence. However, if $k$ is an infinite field, any ideal which is primary to the maximal graded ideal and which is generated by linear forms can be generated by a homogeneous filter-regular sequence (proof of [T1, Lemma 3.1]). \sk

For $i = 0,\ldots,t$ we put
$$a_{\bf z}^i(S/I) := \sup\{r|\ [(I,z_1,\ldots,z_i):z_{i+1}]_r \neq (I,z_1,\ldots,z_i)_r\},$$
with $a_{\bf z}^i(S/I) = -\infty$ if $(I,z_1,\ldots,z_i):z_{i+1} = (I,z_1,\ldots,z_i)$. These invariants can be $\infty$ and they are a measure for how far $\bf z$ is from being a regular sequence in $S/I$. It can be shown that $\bf z$ is a filter-regular sequence for $S/I$ if and only if $a_{\bf z}^i(S/I) < \infty$ for $i = 0,\ldots,t$ [T1, Lemma 2.1]. Note that our definition of $a_{\bf z}^i(S/I)$ is one less than that in [T1]. There is the following close relationship between these numbers and the partial regularity of $S/I$. 

\begin{Theorem} \label{reg1} {\rm [T1, Proposition 2.2]} Let ${\bf z}$ be a filter-regular sequence  of linear forms for $S/I$. Then
$$\reg_t(S/I) = \max\{a_{\bf z}^i(S/I)|\ i = 0,\ldots,t\}.$$ \end{Theorem}\sk

We will use the following characterization of $a_{\bf z}^i(S/I)$.  

\begin{Lemma} \label{other1} $a_{\bf z}^i(S/I) = \max\{r|\ [\cup_{m\ge 1}(I,z_1,\ldots,z_i):z_{i+1}^m]_r \neq (I,z_1,\ldots,z_i)_r\}.$  \end{Lemma}

\begin{pf} Put $r_0 = \max\{r|\ [\cup_{m\ge 1}(I,z_1,\ldots,z_i):z_{i+1}^m]_r \neq (I,z_1,\ldots,z_i)_r\}.$ By definition,
$a_{\bf z}^i(S/I) \le r_0$. Conversely, if $y$ is an element of $\cup_{m\ge 1} (I,z_1,\ldots,z_i):z_{i+1}^m]_{r_0}$, then 
$$yz_i \in [\cup_{m\ge 1}(I,z_1,\ldots,z_i):z_{i+1}^m]_{r_0+1} =   (I,z_1,\ldots,z_i)_{r_0+1}.$$
Hence $y \in [(I,z_1,\ldots,z_i):z_{i+1}]_{r_0}$. This implies $r_0 \le a_{\bf z}^i(S/I)$. So we get $r_0 = a_{\bf z}^i(S/I)$.   \end{pf}

Since $\reg(S/I) = \reg_d(S/I)$, to compute $\reg(S/I)$ we need a filter-regular sequence of linear forms of length $d+1$. But that can be avoided by the following observation. 

\begin{Lemma} \label{parameter} Let ${\bf z} = z_1,\ldots,z_d$ be a filter-regular sequence for $S/I$, $d = \dim(S/I)$. Then $\bf z$ is a system of parameters for $S/I$. \end{Lemma}

\begin{pf} Let $\frak p$ be an arbitrary  associated prime ${\frak p}$ of $(I,z_1,\ldots,z_i)$ with $\dim S/{\frak p} = d-i$, $i = 0,\ldots,d-1$. Then ${\frak p} \neq \frak m$ because $\dim S/{\frak p} > 0$. By the definition of a filter-regular sequence, $z_{i+1} \not\in {\frak p}$. Hence $\bf z$ is a homogeneous system of parameters for $S/I$. \end{pf} 

If $\bf z$ is a homogeneous system of parameters for $S/I$, then $S/(I,z_1,\ldots,z_d)$ is of finite length. Hence $ (S/(I,z_1,\ldots,z_d))_r = 0$ for $r$ large enough. Following [NR] we call
$$r_{\bf z}(S/I) := \max\{r|\ (S/(I,z_1,\ldots,z_d))_r \neq 0\}$$
the {\it reduction number} of $S/I$ with respect to $\bf z$.  It is equal to the maximum degree of the generators of $S/I$ as a module over $k[z_1,\ldots,z_d]$ [V]. Note that the minimum of $r_{\bf z}(S/I)$ is called the reduction number of $S/I$. 

\begin{Theorem} \label{reg2} {\rm [BS, Theorem 1.10], [T1, Corollary 3.3]} Let ${\bf z}$ be a filter-regular sequence of $d$ linear forms for $S/I$. Then
$$\reg(S/I) = \max\{a_{\bf z}^0(S/I),\ldots,a_{\bf z}^{d-1}(S/I),r_{\bf z}(S/I)\}.$$ \end{Theorem}\sk

\noindent{\it Remark.} Theorem \ref{reg2} was proved  in [BS] under an additional condition on the maximum degree of the generators of $I$.\sk 

\section{Evaluations of the initial ideal}

Let $c_i(I)$, $i = 0,\ldots,d$, and $r(I)$ be the invariants defined in Section 1 by means of evaluations of $\In(I)$, where $\In(I)$ is the initial ideal of $I$ with respect to the reverse lexicographic order. We will use the results of Section 2 to express $\reg_t(S/I)$ and $\reg(S/I)$ in terms of $c_i(I)$ and $r(I)$. 

\begin{Lemma} \label{other2} For  ${\bf z} = x_n,\ldots,x_{n-t}$ and $i = 0,\ldots,t$ we have
$$a_{\bf z}^i(S/I) = c_i(I).$$ \end{Lemma}

\begin{pf} By [BS, Lemma (2.2)], $[(I,x_n,\ldots,x_{n-i+1}):x_{n-i}]_r = (I,x_n,\ldots,x_{n-i+1})_r$ if and only if $[(\In(I),x_n,\ldots,x_{n-i+1}):x_{n-i}]_r = (\In(I),x_n,\ldots,x_{n-i+1})_r$ for all $r \ge 0$. Therefore
$$a_{\bf z}^i(S/I) = a_{\bf z}^i(S/\In(I)).$$
By Lemma \ref{other1} we get
$$a_{\bf z}^i(S/\In(I)) = 
\sup\{r|\ [\cup_{m\ge 1}(\In(I),x_n,\ldots,x_{n-i+1}):x_{n-i}^m]_r \neq (\In(I),x_n,\ldots,x_{n-i+1})_r\}.$$
Note that $J_i$ is the ideal of $S_i = k[x_1,\ldots,x_{n-i}]$ obtained from $\In(I)$ by the evaluation $x_{n-i+1} = \cdots = x_n = 0$ and that this evaluation corresponds to the canonical isomorphism $S/(x_{n-i+1},\ldots,x_n) \cong S_i$. Then we may rewrite the above formula as
$$a_{\bf z}^i(S/\!\In(I)) = \sup\{r|\ [\cup_{m\ge 1}J_i:x_{n-i}^m]_r \neq (J_i)_r\}.$$
Since $J_i$ is a monomial ideal, $\cup_{m\ge 1}J_i:x_{n-i}^m$ is generated by the monomials $g$ in the variables $x_1,\ldots,x_{n-i-1}$ for which there exists an integer $m \ge 1$ such that $gx_{n-i}^m \in J_i$. Such a monomial $g$ is determined by the condition $g \in \tilde J_i$. Hence  
$$a_{\bf z}^i(S/\In(I)) = \sup\{r|\ (\tilde J_i)_r \neq (J_i)_r\} = c_i(I).$$ \end{pf} 

As a consequence of Lemma \ref{other2} we can use the invariants $c_i(I)$ to check when $x_n,\ldots,x_{n-t}$ is a regular resp.~filter-regular sequence for $S/I$. 

\begin{Corollary} \label{divisor} $x_{n-i}$ is a non-zerodivisor in $S/(I,x_n,\ldots,x_{n-i+1})$  if and only if  $c_i(I) = -\infty$. \end{Corollary}

\begin{pf} By definition, $a_{\bf z}^i(S/I) = -\infty$ if $x_{n-i}$ is a non-zerodivisor in $S/(I,x_n,\ldots,x_{n-i+1})$. Hence the conclusion follows from Lemma \ref{other2}. \end{pf}

\begin{Corollary} \label{filter} Let ${\bf z} = x_n,\ldots,x_{n-t}$. Then ${\bf z}$ is a filter-regular sequence for $S/I$ if and only if $c_i(I) < \infty$ for $i = 0,\ldots,t$. \end{Corollary}

\begin{pf} It is known that  $\bf z$  is a filter-regular sequence for $S/I$  if and only if $a_{\bf z}^i(S/I) < \infty$ for $i = 0,\ldots,t$  [T1, Lemma 2.1].  \end{pf}

Now we can characterize $\reg_t(S/I)$ as follows.

\begin{Theorem} \label{main3} Assume that $c_i(I) < \infty$ for $i = 0,\ldots,t$. Then
$$\reg_t(S/I) = \max\{c_i(I)|\ i = 0,\ldots,t\}.$$ \end{Theorem}

\begin{pf} This follows from Theorem \ref{reg1}, Lemma \ref{other2} and Corollary \ref{filter}. \end{pf}

We can also give a characterization of $\reg(S/I)$ which involves $r(I)$.

\begin{Lemma} \label{other3} Assume that $c_i(I) < \infty$ for $i = 0,\ldots,d-1$. Then
$$r_{\bf z}(S/I) = r(I).$$  \end{Lemma}

\begin{pf} By Corollary \ref{filter}, ${\bf z} = x_n,\ldots,x_{n-d+1}$ is a filter-regular sequence for $S/I$. By Lemma \ref{parameter} and [T2, Theorem 4.1], this implies that $\bf z$ is a homogeneous system of parameters for $S/\In(I)$ with
$$r_{\bf z}(S/I) = r_{\bf z}(S/\In(I)).$$
Note that $S/(x_{n-d+1},\ldots,x_n) \cong S_d$ and that $J_d$ is the ideal obtained from $\In(I)$ by the evaluation $x_{n-d+1} = \cdots = x_n = 0$. Then
\begin{eqnarray*} r_{\bf z}(S/\In(I)) & = & \max\{r|\ (S/(\In(I),x_n,\ldots,x_{n-d+1}))_r \neq 0\}\\
& = & \max\{r|\ (S_d/J_d)_r \neq 0\}\\
& = &  r(I). \end{eqnarray*} 
\end{pf} 

\begin{Theorem} \label{main4} Assume that $c_i(I) < \infty$ for $i = 0,\ldots,d-1$. Then
$$\reg(S/I) = \max\{c_0(I),\ldots,c_{d-1}(I),r(I)\}.$$ \end{Theorem}

\begin{pf} This follows from Theorem \ref{reg2}, Lemma \ref{other2}, Corollary \ref{filter} and Lemma \ref{other3}.\end{pf}

\section{Combinatorial description} 

First, we want to show that the condition $c_i(I) < \infty$ can be easily checked in terms of the lattice vectors of the generators of $\In(I)$.
Let $\cal B$ be the (finite) set of monomials which minimally generates $\In(I)$. We set
$$E_i := \{v \in {\Bbb N}^{n-i}|\ x^v \in {\cal B}\},$$ 
where $x^v = x_1^{\varepsilon_1}\cdots x_s^{\varepsilon_s}$ if $v = (\varepsilon_1,\ldots,\varepsilon_s)$. For $j = 1,\ldots,n-i$ we denote by $p_j$ the projection from ${\Bbb N}^{n-i}$ to ${\Bbb N}^{n-i-1}$ which deletes the $j$th coordinate. For two lattice vectors $a = (\alpha_1,\ldots,\alpha_s)$ and $b = (\beta_1,\ldots,\beta_s)$ of the same size we say $a \ge b$ if $\alpha_j \ge \beta_j$ for $j = 1,\ldots,s$. 

\begin{Lemma} \label{finite} $c_i(I) < \infty$ if and only if for every element $a \in p_{n-i}(E_i)\setminus E_{i+1}$ there are elements $b_j \in E_{i+1}$ such that $p_j(a) \ge p_j(b_j)$, $j = 1,\ldots,n-i-1$. \end{Lemma}

\begin{pf} Recall that $c_i(I) = \sup\{r|\ (\tilde J_i/J_i)_r \neq 0\}$. Then $c_i(I) < \infty$ if and only if $\tilde J_i/J_i$ is of finite length. By the definition of $J_i$ and $\tilde J_i$, the latter condition is equivalent to the existence of a number $r$ such that $x_j^r\tilde J_i \subseteq J_i$ for $j = 1,\ldots,n-i$. It is clear that $J_i$ is generated by the monomials $x^v$ with $v \in E_i$. From this it follows that $\tilde J_i$ is generated by $J_i$ and the monomials $x^a$ with $a \in p_{n-i}(E_i)\setminus E_{i+1}$. For such a monomial $x^a$ we can always find a number $r$ such that $x_{n-i}^rx^a \in J_i$. For $j < n-i$, $x_j^rx^a \in J_i$ if and only if  $x_j^rx^a$ is divided by a generator $x^{b_j}$ of $J_i$. Since $x_j^rx^a$ does not contain $x_{n-i},...,x_n$, so does $x^{b_j}$. Hence $b_j \in E_{i+1}$. Setting $x_j=1$ we see that $x_j^rx^a$ is divided by $x^{b_j}$ for some number $r$ if and only if $p_j(a) \ge p_j(b_j)$. \end{pf}

If $c_i(I) = \infty$, we should make a random linear transformation of the variables $x_1,\ldots,x_{n-i}$ and test the condition $c_i(I) < \infty$ again. By Lemma \ref{other2} the linear transformation does not change the invariants $c_j(I)$ for $j < i$. Moreover, instead of $\In(I)$ we only need to compute the smaller initial ideal $\In(I_i)$, where $I_i$ denotes the ideal of $S_i$ obtained from $I$ by the evaluation $x_{n-i+1} = \cdots = x_n = 0$. Let ${\cal B}_i$ be the set of monomials which minimally generates $\In(I_i)$. It is easy to see that ${\cal B}_i$ is the set of the monomials of $\cal B$ which are not divided by $x_{n-i+1},\ldots,x_n$. From this it follows that $E_j = \{v \in {\Bbb N}^{n-j}|\ x^v \in {\cal B}_i\}$ for $j \le i$. Thus, we can use this formula to compute $E_j$ and to check the condition $c_j(I) < \infty$ for $j \le i$. Once we know $c_i(I) < \infty$ we can proceed to compute $c_i(I)$. \sk

In the lattice ${\Bbb N}^{n-i}$ we delete the shadow of $E_i$, that is, the set of elements $a$ for which there is $v \in E_i$ with $v \le a$. The remaining lattice has the shape of a staircase and we will denote by $F_i$ the set of its corners. It is easy to see that $F_i$ is the set of the elements of the form $a = \max(v_1,\ldots,v_{n-i})-(1,\ldots,1)$ 
with $a \not\ge v$ for any element $v \in E_i$, where $v_1,\ldots,v_{n-i}$ is a family of $n-i$ elements of $E_i$ for which the $j$th coordinate of $v_j$ is greater than the $j$th coordinate of $v_h$ for all $h \neq j$, $j = 1,\ldots,n-i$, and $\max(v_1,\ldots,v_{n-i})$ denotes the element whose coordinates are the maxima of the corresponding coordinates of $v_1,\ldots,v_{n-i}$. If $a = (\alpha_1,\ldots,\alpha_{n-i})$, we set
$$|a| := \alpha_1+\ldots+\alpha_{n-i}.$$

\begin{Proposition} \label{other4}  Assume that $c_i(I) < \infty$. Then $c_i(I) = -\infty$ if $F_i = \emptyset$ and $c_i(I) =  \max_{a \in F_i}|a|$ if $F_i \neq -\emptyset$ . \end{Proposition}

\begin{pf} Let $a$ be an arbitrary element of $F_i$. Then $a = \max(v_1,\ldots,v_{n-i})-(1,\ldots,1)$ for some family  $v_1,\ldots,v_{n-i}$ of $S_i$. Let $v_j = (\varepsilon_{j1},\ldots,\varepsilon_{jn-i})$, $j = 1,\ldots,n-i$. Then $a = (\varepsilon_{11}-1,\ldots,\varepsilon_{n-in-i}-1)$. Since $\varepsilon_{jj} > \varepsilon_{hj}$ for $h \neq j$, we get 
$a \ge (\varepsilon_{n-i1},\ldots,\varepsilon_{n-in-i-1},0).$ Therefore, $x^a$ is divided by the monomial obtained from $x^{v_{n-i}}$ by setting $x_{n-i} = 1$. Note that $J_i$ is generated by the monomials $x^v$ with $x_v \in E_i$. Since $v_{n-i} \in E_i$, we have $x^{v_{n-i}} \in J_i$, whence $x^a \in \tilde J_i$. On the other hand, $x^a \not\in J_i$ because $a \not\ge v$ for any element $v \in E_i$. Since $|a| = \deg x^a$, this implies $(\tilde J_i/J_i)_{|a|} \neq 0$. Hence $|a| \leq c_i(I)$. So we obtain $\max_{a\in F_i}|a| \le c_i(I)$ if $F_i \neq \emptyset$.  \par
To prove the converse inequality we assume that $\tilde J_i/J_i \neq 0$. Since $c_i(I) < \infty$, there is a monomial $x^b \in \tilde J_i \setminus J_i$ such that $\deg x^b = c_i(I)$. Since $x^b \not\in J_i$, $b \not\ge v$ for any element $v \in E_i$. For $j = 1,\ldots,n-i$ we have $x_jx^b \in J_i$ because $\deg x_jx^b = c_i(I)+1$. Therefore, $x_jx^b$ is divided by some monomial $x^{v_j}$ with $v_j \in E_i$. Let $b = (\beta_1,\ldots,\beta_{n-i})$ and $v_j = (\varepsilon_{j1},\ldots,\varepsilon_{jn-i})$. Then $\beta_h \ge \varepsilon_{jh}$ for $h \neq j$ and $\beta_j + 1 \ge \varepsilon_{jj}$. 
Since $b \not\ge v_j$, we must have $\beta_j < \varepsilon_{jj}$, hence $\beta_j = \varepsilon_{jj}-1$. It follows that $\varepsilon_{jj} = \beta_j + 1 > \varepsilon_{hj}$ for all $h \neq j$. Thus, the family $v_1,\ldots,v_{n-i}$ belongs to ${\cal S}_i$ and $b = \max(v_1,\ldots,v_{n-i})-(1,\ldots,1)$. So we have proved that $b \in F_i$.  Hence $c_i(I) = \deg x^b = |b|  \le \max_{a \in F_i}|a|.$ \par
The above argument also shows that $F_i \neq \emptyset$ if $\tilde J_i \neq J_i$. So $c_i(I) = -\infty$ if $F_i = \emptyset$.\end{pf}
 
By Corollary \ref{filter}, if $c_i(I) < \infty$ for $i = 0,\ldots,d-1$, then ${\bf z} = x_n,\ldots,x_{n-d+1}$ is a filter-regular sequence for $S/I$. By Lemma \ref{parameter} and Lemma \ref{other3}, that implies $r(I) = r_{\bf z}(S/I) < \infty$. In this case, we have the following description of $r(I)$.

\begin{Proposition} \label{other5}  Assume that $r(I) < \infty$. Then $r(I) = \max_{a \in F_d}|a|$. \end{Proposition}

\begin{pf} This can be proved similarly as in the proof of Lemma \ref{other4}. \end{pf}

Combining the above results with Theorem \ref{main3} and Theorem \ref{main4} we get a simple method to compute $\reg_t(S/I)$ and $\reg(S/I)$. We will illustrate the above method by an example at the end of the next section. Moreover, we get the following estimation for $\reg_t(S/I)$. 

\begin{Corollary} \label{bound} Let $x_n,\ldots,x_{n-t}$ be a filter-regular sequence for $S/I$. Let $g_i$ denote the least common multiple of the minimal generators of $\In(I)$ which are not divided by any of the variables $x_{n-i+1},\ldots,x_n$. Then
$$\reg_t(S/I) \le \max\{\deg g_i-n+i|\ i=0,\ldots,t\}.$$\end{Corollary}

\begin{pf} By Corollary \ref{filter}, the assumption implies that $c_i(I) < \infty$ for $i = 0,...,t$. Thus, combining Theorem \ref{main3} and Lemma \ref{other4} we get
$$\reg_t(S/I) \le \max\{|a||\ a \in F_i,\ i = 0,\ldots,t\}.$$
It is easily seen from the definition of $F_i$ that $\max_{a\in F_i}|a| \le \deg g_i-n+i$, $i = 0,\ldots,t$, hence the conclusion. \end{pf}

\noindent{\it Remark.} Bruns and Herzog [BH, Theorem 3.1(a)] resp.~Hoa and Trung [HT, Theorem 3.1] proved that for any monomial ideal $I$, $\reg(S/I) \le \deg f-1$ resp.~$\deg f - \height I$, where $f$ is the least common multiple of the minimal generators of $I$. Note that the mentioned result of Bruns and Herzog is valid for multigraded modules. \sk
 
\section{The case of projective curves}

Let $I_C \subset k[x_1,\ldots,x_n]$ be the defining  saturated ideal of a (not necessarily reduced) projective curve $C \subset {\Bbb P}^{n-1}$, $n \ge 3$. We will assume that $k[x_{n-1},x_n] \hookrightarrow S/I_C$ is a Noether normalization of $S/I_C$. In this case, Theorem \ref{main4} can be reformulated as follows.

\begin{Proposition} \label{curve1}  $\reg(S/I_C) = \max\{c_1(I_C),r(I_C)\}.$  \end{Proposition}

\begin{pf} By the above assumption $S/I_C$ is a generalized Cohen-Macaulay ring of positive depth and $x_n,x_{n-1}$ is a homogeneous system of parameters for $S/I_C$. Therefore, $x_n,x_{n-1}$ is a filter-regular sequence for $S/I_C$. In particular, $x_n$ is a non-zerodivisor in $S/I_C$. By  Lemma \ref{divisor}, $c_0(I_C) = -\infty$. Hence the conclusion follows from Theorem \ref{main4}. \end{pf}

Since $S/I_C$ has positive depth, the graded minimal  free resolution of $S/I_C$ ends at most at the $(n-1)$th place:
$$0 \To F_{n-1} \To \cdots \To F_1 \To F_0 \To S/I_C \To 0.$$
>From Theorem \ref{main3} we obtain the following information on the shifts of $F_{n-1}$. Note that  $F_{n-1} = 0$ if $S/I_C$ is a Cohen-Macaulay ring or, in other words, if $C$ is an arithmetically Cohen-Macaulay curve. 

\begin{Proposition} \label{curve2} If $C$ is not an arithmetically Cohen-Macaulay curve, $c_1(I_C)+n-1$ is the maximum degree of the generators of  $F_{n-1}$. \end{Proposition}

\begin{pf} Let $b_{n-1}$ be the maximum degree of the generators of  $F_{n-1}$. As we have seen in the introduction, 
$b_{n-1}-n+1 = (n-1)\text{--}\!\reg(S/I_C)  = \reg_1(S/I_C).$ By Theorem \ref{main3}, $\reg_1(S/I_C) = \max\{c_0(I_C),c_1(I_C)\} = c_1(I_C)$ because $c_0(I_C) = -\infty$. So we obtain $b_{n-1} = c_1(I_C)+n-1$. \end{pf}

Now we shall see that Proposition \ref{curve1} contains all main results of Bermejo and Gimenez in [BG]. It should be noted that they did not use strong results such as Theorem \ref{reg2}. We follow the notations of [BG].\sk

Let $E := \{a \in {\Bbb N}^{n-2}|\ x^a \in \In(I_C)\}$ and denote by $H(E)$ the smallest integer $r$ such that $a \in E$ if $|a| = r$. 

\begin{Corollary} \label{BG1} {\rm [BG, Theorem 2.4]} Assume that $C$ is an arithmetically Cohen-Macaulay curve. Then $\reg(S/I_C) = H(E)-1.$ \end{Corollary}

\begin{pf} Since $x_n,x_{n-1}$ is a regular sequence in $S/I_C$, we have $c_1(I_C) = -\infty$ by Corollary \ref{divisor}.  By Proposition \ref{curve1} this implies $\reg(S/I_C) = r(I_C)$. But  
$$r(I_C) = \sup\{r|\ (S_2/J_2)_r \neq 0\} =H(E)-1$$
because $J_2$ is generated by the monomials $x^a$, $a \in E$. \end{pf}  

Let $I_0$ be the ideal in $S$ generated by the polynomials obtained from $I_C$ by the evaluation $x_{n-1} = x_n = 0$. Then $S/I_0$ is a two-dimensional Cohen-Macaulay ring. Let $\tilde I$ denote the ideal in $S$ generated by the monomials obtained from $\In(I_C)$ by the evaluation $x_{n-1} = x_n = 1$. Let
$$F := \{a \in{\Bbb N}^{n-2}|\ x^a \in \tilde I \setminus \In(I_0)\}.$$
For every vector $a \in F$ let 
$$E_a := \{(\mu,\nu)\in {\Bbb N}^2|\ x^a x_{n-1}^{\mu}x_n^{\nu} \in \In(I_C)\}.$$
Let $\Re := \cup_{a \in F}\{a\times[{\Bbb N}^2 \setminus E_a]\}$ and denote by $H(\Re)$ the smallest integer $r$ such that the number of the elements $b \in \Re$ with $|b| = s$ becomes a constant for $s \ge r$. 

\begin{Corollary} \label{BG2} {\rm [BG, Theorem 2.7]}  $\reg(S/I_C) = \max\{\reg(S/I_0),H(\Re)\}.$ \end{Corollary}

\begin{pf}  As in the proof of Corollary \ref{BG1} we have  $\reg(S/I_0) = r(I_0)$. But $r(I_0) = r(I_C)$ because $\In(I_0)$ is  the ideal generated by the monomials obtained from $\In(I_C)$ by the evaluation $x_{n-1} = x_n = 0$. Thus, 
$$\reg(S/I_0) = r(I_C).$$
It has been observed in [BG] that the number of the elements $b \in \Re$ with $|b| = s$ is the difference 
$H_{S/I_C}(s)- H_{S/\tilde I}(s) = H_{S/\In(I_C)}(s) - H_{S/\tilde I}(s) = H_{\tilde I/\In(I_C)}(s),$
where $H_E(s)$ denotes the Hilbert function of a graded $S$-module $E$. Since $x_n$ is a non-zerodivisor in $S/\In(I_C)$, $H(\Re)+1$ is the least integer $r$ such that $H_{(\tilde I,x_n)/(\In(I_C),x_n)}(s) = 0$ for $s \ge r$. On the other hand, since $\In(I_C)$ is generated by monomials which do not contain $x_n$ and since $J_1$ is the ideal in $k[x_1,\ldots,x_{n-1}]$ obtained from $\In(I_C)$ by the evaluation $x_n=0$, we have $\In(I_C) = J_1S$ and $\tilde I = \tilde J_1S$, whence $(\tilde I,x_n)/(\In(I_C),x_n) \cong \tilde J_1/J_1$. Note that $c_1(I_C) = \max\{r|\ (\tilde J_1/J_1)_r \neq 0\}$ with $c_1(I_C) = -\infty$ if $\tilde J_1 = J_1$. Then
$$H(\Re) = \max\{0,c_1(I_C)\}.$$
Thus, applying Proposition \ref{curve1} we obtain  $\reg(S/I_C) = \max\{\reg(S/I_0),H(\Re)\}.$
\end{pf}

\noindent{\bf Example.} Let $C \subset {\Bbb P}^3$ be the monomial curve  
$(t^\alpha s^\beta: t^\beta s^\alpha:s^{\alpha+\beta}:t^{\alpha+\beta})$,   $\alpha > \beta > 0$,   g.c.d$(\alpha,\beta) = 1$. It is known that the defining ideal $I_C \subset k[x_1,x_2,x_3,x_4]$ is generated by the quadric $x_1x_2-x_3x_4$ and the forms $x_1^{\beta+r}x_3^{\alpha-\beta-r}-x_2^{\alpha-r}x_4^r$, $r = 0,\ldots,\alpha-\beta$, and that this is a Gr\"obner basis of $I_C$ for  the reverse lexicographic order with $x_1 > x_2 > x_3 > x_4$ [CM, Th\'eor\`em 3.9]. Therefore,
$$\In(I_C) = (x_1x_2,x_2^\alpha,x_1^{\beta+1}x_3^{\alpha-\beta-1},x_1^{\beta+2}x_3^{\alpha-\beta-2},\ldots,x_1^\alpha).$$ Using the notations of Section 3 we have
\begin{eqnarray*} E_1 & = & \{(1,1,0),(0,\alpha,0),(\beta+1,0,\alpha-\beta-1),(\beta+2,0,\alpha-\beta-2),\ldots,(\alpha,0,0)\},\\
E_2 & = & \{(1,1),(0,\alpha), (\alpha,0)\}. \end{eqnarray*}
>From this it follows that
\begin{eqnarray*} F_1 & = & \{(\beta+1,0,\alpha-\beta-2), (\beta+2,0,\alpha-\beta-3),\ldots,(\alpha-1,0,0)\},\\
F_2 & = & \{(0,\alpha-1),(\alpha-1,0)\}.\end{eqnarray*}
By Proposition \ref{other4}, $c_1(I_C) = \alpha-1$ if $\alpha -\beta \ge 2$ ($c_1(I_C) = -\infty$ if $\alpha -\beta = 1$) and $r(I_C) = \alpha-1$ by Proposition \ref{other5}. Applying Proposition \ref{curve1} we obtain $\reg(S/I_C) = \alpha-1$.\par

The direct computation of the invariant $H(\Re)$ is more complicated than that of $c_1(I_C)$. First, we should interpret $F$ as the set of the elements of the form $a \in {\Bbb N}^2$ such that $a \ge b$ for some elements $b \in p(E_1)$ but $a \not\ge c$ for any element $c \in E_2$. Then we get
$$F = \{(\beta+1,0),(\beta+2,0),\ldots,(\alpha-1,0)\}.$$
For all $\varepsilon = \beta+1,\ldots,\alpha-1$ we verify that $E_{(\varepsilon,0)} = (\alpha-\varepsilon,0)+{\Bbb N}^2$. It follows that 
$$\Re = \{(\varepsilon,0,\mu,\nu)\in {\Bbb N}^4|\ \varepsilon = \beta+1,\ldots,\alpha-1;\, \mu \le \alpha-\varepsilon-1\}.$$
If $\alpha-\beta = 1$, we have $\Re = \emptyset$, hence $H(\Re) = 0$. If $\alpha-\beta \ge 2$, we can check that $H(\Re) =  \alpha-1.$  \sk

\section*{References}

\noindent [AH] A. Aramova and J. Herzog, Almost regular sequences and Betti numbers, preprint.\par
\noindent [BCP] D.~Bayer, H.~Charalambous and S. Popescu, Extremal Betti numbers and applications to monomial ideals, J. Algebra 221 (1999), 497-512. \par
\noindent [BM] D.~Bayer and D.~Mumford, What can be computed in algebraic geometry~? in: D. Eisenbud and L. Robbiano (eds.), Computational Algebraic
Geometry and Commutative Algebra, Proceedings, Cortona (1991), Cambridge
University Press, 1993, 1-48.\par
\noindent [BS] D.~Bayer and M.~Stillman, A criterion for detecting $m$-regularity, Invent. Math. 87 (1987), 1-11.\par
\noindent [BG] I.~Bermejo and P. Gimenez, On the Castelnuovo-Mumford regularity of projective curves, Proc. Amer. Math. Soc. 128 (2000), 1293-1299. \par
\noindent [BH] W.~Bruns and J.~Herzog, On multigraded resolutions, Math. Proc. Cambridge Phil. Soc. 118 (1995), 245-275.\par
\noindent [CM] L.~Coudurier and M.~Morales, Classification des courbes toriques dans l'espace projectif, module de Rao et liaison, J.~Algebra 211 (1999), 524-548.\par
\noindent [EG] D.~Eisenbud and S.~Goto, Linear free resolutions and minimal multiplicities, J.~Algebra~88 (1984), 89-133.\par
\noindent [HT] L.T.~Hoa and N.V.~Trung, On the Castelnuovo-Mumford regularity and the arithmetic degree of monomial ideals, Math.~Z. 229 (1998), 519-537.\par
\noindent [NR] D.G.~Northcott and D.~Rees, Reductions of ideals in local rings, Proc. Cambridge Phil. Soc. 50 (1954), 145-158. \par
\noindent [STC]  P.~Schenzel, N.V.~Trung and N.T.~Cuong, \"Uber verallgemeinerte Cohen-Macaulay-Moduln,  Math. Nachr. 85 (1978), 57-73.\par
\noindent [T1] N.V.~Trung, Reduction exponent and degree bounds for the defining equations of a graded ring, Proc. Amer. Math. Soc. 102 (1987), 229-236.\par
\noindent [T2] N.V.~Trung, Gr\"obner bases, local cohomology and reduction number, to appear in Proc. Amer. Math. Soc.\par
\noindent [V] W.~Vasconcelos, Cohomological degree of a module, in: J.~Elias, J.M.~Giral, R.M.~Miro-Roig, S.~Zarzuela (eds.), Six Lectures on Commutative Algebra, Progress in Mathematics 166, pp.~345-392, Birkh\"auser, 1998.
\end{document}